\RequirePackage[british]{babel}
\documentclass[reqno,a4paper,12pt]{amsart}
\usepackage{a4wide}
\usepackage[utf8x]{inputenc}
\usepackage{amsmath,amssymb,amstext,amsthm,amscd,mathrsfs,eucal}
\usepackage{graphicx,color}
\usepackage{array}
\usepackage[noadjust]{cite}
\usepackage{hyperref}
\hypersetup{%
  pdftitle   = {Lorentzian Lie n-algebras},
  pdfkeywords = {lorentzian, M2, M Theory, Bagger-Lambert, Lie, Filippov, k-algebra},
  pdfauthor  = {José Figueroa-O'Farrill},
  pdfcreator = {\LaTeX\ with package \flqq hyperref\frqq}
}
\PrerenderUnicode{éÉ}
\newcommand{\half}{\tfrac12}

\newcommand{\fa}{\mathfrak{a}}
\newcommand{\fb}{\mathfrak{b}}

\newcommand{\fs}{\mathfrak{s}}
\newcommand{\fsl}{\mathfrak{sl}}

\newcommand{\RR}{\mathbb{R}}

\newcommand{\be}{\boldsymbol{e}}
\newcommand{\bx}{{\boldsymbol{x}}}

\DeclareMathOperator{\End}{End}

\DeclareMathOperator{\ad}{ad}

\theoremstyle{plain}
\newtheorem{lemma}{Lemma}
\newtheorem{proposition}[lemma]{Proposition}
\newtheorem{theorem}[lemma]{Theorem}
\newtheorem{corollary}[lemma]{Corollary}
\theoremstyle{definition}

\newcommand{\MUNCH}[1]{\relax}

\allowdisplaybreaks
\begin{document}
\title{Lorentzian Lie $n$-algebras}
\author{José Figueroa-O'Farrill}
\address{Maxwell Institute and School of Mathematics, University of
  Edinburgh, UK}
\email{J.M.Figueroa@ed.ac.uk}
\date{\today}
\begin{abstract}
  We classify Lie $n$-algebras possessing an invariant lorentzian
  inner product.
\end{abstract}
\maketitle
\tableofcontents

\section{Introduction}
\label{sec:intro}

It is natural to classify algebraic structures, especially those which
do not arise in Algebra alone.  Such is the case with Lie $n$-algebras
(see below for a definition), which have made their appearance in
several physical contexts independently from the initial purely
algebraic aim of generalising the notion of a Lie algebra.

General classifications, however, are usually hard, unless one imposes
further conditions.  For example, the classification problem of Lie
algebras is \emph{not tame}, which in layman's terms translates as
being hopeless.  Standing on the shoulders of Killing and Élie Cartan,
the semisimple Lie algebras can be classified in a relatively
straightforward manner, whereas the Levi decomposition further reduces
the problem to classifying solvable Lie algebras and their possible
semidirect products with semisimple Lie algebras.  However the
classification of solvable Lie algebras except in very low dimension
is not possible.  Even nilpotent Lie algebras stop being tame in
dimension $\geq 7$.

One natural class of Lie algebras which extend the semisimple Lie
algebras and which, due to their relative rarity, offer a hope of
classification, are the metric Lie algebras: those possessing a
(nondegenerate) ad-invariant inner product.  It is a classical result
that the metric Lie algebras for which the inner product is
positive-definite are the reductive Lie algebras, which are generated
under direct sum by the one-dimensional and the simple Lie algebras.
For indefinite signature, things get a little more interesting.
Medina and Revoy \cite{MedinaRevoy} (see also \cite{FSalgebra}) showed
that metric Lie algebras of any signature are still generated by the
same ingredients under direct sum and a new construction called a
\emph{double extension}.  This construction has its origin in the
classification of lorentzian Lie algebras, summarised as follows.

\begin{theorem}[\cite{MedinaLorentzian}]\label{th:lorentzianLie}
  A finite-dimensional indecomposable lorentzian Lie algebra
  is either one-dimensional, isomorphic to $\fsl(2,\RR)$ or else
  isomorphic to $\RR u \oplus \RR v \oplus W$, with Lie bracket
  \begin{equation*}
    [u,x] = A(x) \qquad\text{and}\qquad [x,y] = \left<A(x), y\right> v~,
  \end{equation*}
  where $\left<-,-\right>$ is a positive-definite inner product on
  $W$, $A: W \to W$ is an invertible skewsymmetric endomorphism and we
  extend the inner product on $W$ to all of $V$ by declaring $u,v
  \perp W$, $\left<u,v\right> =1$ and $\left<u,u\right> =
  \left<v,v\right> = 0$.
\end{theorem}

Pushing this construction further allows a construction of all metric
Lie algebras with signature $(2,p)$ by double extending the lorentzian
Lie algebras \cite{MedinaRevoy}; however, the description of such
metric Lie algebras as a double extension is ambiguous and it is not
clear a priori which of the Lie algebras so obtained are isomorphic.
A different approach which results in a classification of Lie algebras
with signature $(2,p)$ is given by Kath and Olbrich in
\cite{KathOlbrich2p}, following the announcement \cite{baumkath2p}.

Lie $n$-algebras (for $n>2$) are a natural generalisation of the
notion of a Lie algebra, with which they coincide when $n=2$.  They
arise in a number of unrelated problems in mathematical physics in
areas such as integrable systems, supersymmetric gauge theories,
supergravity and string theory.  It is this which, to my mind, makes
their classification into a \emph{de facto} interesting problem.

A \textbf{Lie $n$-algebra} structure on a vector space $V$ is a linear
map $\Phi: \Lambda^n V \to V$, denoted simply by a $n$-ary bracket,
with the property that for every $x_1,\dots,x_{n-1} \in V$, the left
multiplication $\ad_{x_1,\dots,x_{n-1}}: y \mapsto
[x_1,\dots,x_{n-1},y]$ is a derivation over the bracket.  If $n=2$
this latter property is precisely the Jacobi identity for a Lie
algebra.  For $n>2$ we will call it the $n$-Jacobi identity.  For
$n>2$, Lie $n$-algebras were introduced by Filippov \cite{Filippov}
and are often referred to as Filippov $n$-algebras.  Many of the
structural results in the theory of Lie algebras have their analogue
in the theory of Lie $n$-algebras, often with some refinement;
although it seems that Lie $n$-algebras become more and more rare as
$n$ increases, due perhaps to the fact that the $n$-Jacobi identity
imposes more and more conditions as $n$ increases.  From now on,
whenever we write Lie $n$-algebra, we will assume that $n>2$ unless
otherwise stated.

For example, over the complex numbers there is up to isomorphism a
unique simple Lie $n$-algebra for every $n>2$, of dimension $n+1$ and
whose $n$-bracket is given relative to a basis $(\be_i)$ by
\begin{equation*}
  [\be_1,\dots,\widehat{\be_i},\dots,\be_{n+1}] = (-1)^i \be_i~,
\end{equation*}
where a hat denotes omission.  Over the reals, they are all given by
attaching a sign $\varepsilon_i$ to each $\be_i$ on the right-hand
side of the bracket.  This result is due to Ling \cite{LingSimple},
who also established a Levi decomposition of an arbitrary Lie
$n$-algebra into a direct sum of a semisimple Lie $n$-algebra (a
direct sum of its simple ideals) and a maximal solvable subalgebra.
This again reduces the classification of Lie $n$-algebras to that of
the solvable Lie $n$-algebras and their semidirect products with
semisimple Lie $n$-algebras, but classifying solvable Lie $n$-algebras
seems to be as hard for $n>2$ as it is for $n=2$.

Taking a cue from the case of Lie algebras, and because it is this
class of Lie $n$-algebras which seem to appear in nature, one can
restrict to metric Lie $n$-algebras.  Let $b \in S^2V^*$ be an inner
product (i.e., a nondegenerate symmetric bilinear form), denoted
simply as $\left<-,-\right>$.  We say that a Lie $n$-algebra
$(V,\Phi,b)$ is \textbf{metric} if the left multiplications
$\ad_{x_1,\dots,x_{n-1}}$ for $x_i\in V$ are skewsymmetric relative to
$b$.  Metric Lie $n$-algebras seem to have appeared for the first time
in \cite{FOPPluecker}, albeit tangentially, and more prominently in
\cite{BL2}, which has attracted a great deal of attention recently in
the mathematical/theoretical physics community.

Given two metric Lie $n$-algebras $(V_1,\Phi_1,b_1)$ and
$(V_2,\Phi_2,b_2)$, we may form their \textbf{orthogonal direct sum}
$(V_1\oplus V_2,\Phi_1\oplus \Phi_2, b_1 \oplus b_2)$, by declaring
that
\begin{align*}
  [x_1,x_2,y_1,\dots,y_{n-2}] = 0 \qquad\text{and}\qquad
  \left<x_1,x_2\right> = 0~,
\end{align*}
for all $x_i\in V_i$ and all $y_i\in V_1 \oplus V_2$.  The resulting
object is again a metric Lie $n$-algebra.  A metric Lie $n$-algebra is
said to be indecomposable if it is not isomorphic (see below) to an
orthogonal direct sum of metric Lie $n$-algebras $(V_1\oplus 
V_2, \Phi_1\oplus \Phi_2, b_1\oplus b_2)$ with $\dim V_i > 0$.
In order to classify the metric Lie $n$-algebras, it is clearly enough
to classify the indecomposable ones.

The orthogonal Plücker identities conjectured in \cite{FOPPluecker}
imply that indecomposable Lie $n$-algebras admitting a
positive-definite inner product are either simple or one-dimensional.
This conjecture was proved by Nagy in \cite{NagykLie} and
independently by Papadopoulos \cite{GPkLie}.  For the special case of
$n=3$ this result was rediscovered in \cite{GP3Lie,GG3Lie}.
Lorentzian Lie $3$-algebras have been classified in \cite{Lor3Lie},
where a one-to-one correspondence is established between
indecomposable lorentzian Lie $3$-algebras and semisimple euclidean
Lie algebras.  More precisely one has the following

\begin{theorem}[\cite{Lor3Lie}]\label{th:lorentzian3Lie}
  Let $(V,\Phi,b)$ be a finite-dimensional indecomposable lorentzian
  Lie 3-algebra.  Then it is either one-dimensional, simple or else
  isomorphic to $\RR u \oplus \RR v \oplus W$ with 3-bracket
  \begin{equation*}
    [u,x,y] = [x,y] \qquad\text{and}\qquad [x, y, z] = -
    \left<[x,y],z\right> v~,
  \end{equation*}
  where $[-,-]: \Lambda^2 W \to W$ is a Lie bracket making $W$ into a
  compact semisimple Lie algebra and $\left<-,-\right>$ is a
  positive-definite ad-invariant inner product, extended to all of $V$
  by declaring $u,v \perp W$, $\left<u,v\right> =1$ and
  $\left<u,u\right> = \left<v,v\right>= 0$.
\end{theorem}

This latter class of Lie 3-algebras were discovered independently in
\cite{GMRBL,BRGTV,HIM-M2toD2rev}.

The purpose of this short note is to classify lorentzian Lie
$n$-algebras for $n>3$.   The classification borrows from the
techniques in \cite{Lor3Lie}, which themselves are inspired by
\cite{MedinaRevoy,FSalgebra}.  The main result is contained in
Theorem \ref{th:lorentzianNLie}.  The classification of metric Lie
$n$-algebras will not be attempted here, but see \cite{2p3Lie} for the
case $n=3$ and signature $(2,p)$.

\section*{Acknowledgments}

It is a pleasure to thank Paul de Medeiros and Elena Méndez-Escobar
for many $n$-algebraic discussions.  I would also like to thank the
combined efforts of Martin Frick, Wolfgang Hein and Nils-Peter
Skoruppa of the Fachbereich Mathematik der Universität Siegen for
making available a copy of \cite{LingSimple}.

\section{Metric Lie $n$-algebras}
\label{sec:metric-n-lie}

We will require the basic terminology of Lie $n$-algebras, as in
\cite{Filippov} or \cite{LingSimple} and of metric Lie $n$-algebras.
Many of the proofs will be omitted, since they can be read
\emph{mutatis mutandis} from the ones for $n=3$ in \cite{Lor3Lie}.

\subsection{Some structure theory}
\label{sec:structure}

Let $(V,\Phi)$ be a Lie $n$-algebra.  Given subspaces $W_i \subset V$,
we will let $[W_1,\dots, W_n]$ denote the subspace of $V$ consisting
of elements $[w_1,\dots,w_n] \in V$, where $w_i \in W_i$.

A subspace $W \subset V$ is a \textbf{subalgebra}, written $W < V$, if
$[W,\dots,W]\subset W$.  A subalgebra $W<V$ is said to be
\textbf{abelian} if $[W,\dots,W]=0$.

If $V,W$ are Lie $n$-algebras, then a linear map $\phi: V \to W$ is a
\textbf{homomorphism} if
\begin{equation*}
  \phi[x_1,\dots,x_n] = [\phi(x_1), \dots, \phi(x_n)]~,
\end{equation*}
for all $x_1,\dots,x_n\in V$.  If $\phi$ is also a vector space
isomorphism, we say that it is an \textbf{isomorphism} of Lie
$n$-algebras.  An isomorphism $V \to V$ is called an
\textbf{automorphism}.  A subspace $I \subset V$ is an \textbf{ideal},
written $I \lhd V$, if $[I,V,\dots,V]\subset I$.  It follows that
there is a one-to-one correspondence between ideals and kernels of
homomorphisms.  If $I \lhd V$ and $J\lhd V$, then $I\cap J \lhd V$ and
$I + J \lhd V$.  We will say that an ideal $I \lhd V$ is
\textbf{minimal} if any other ideal $J \lhd V$ contained in $I$ is
either $0$ or $I$.  Dually, an ideal $I \lhd V$ is \textbf{maximal} if
any other ideal $J \lhd V$ containing $I$ is either $I$ or $V$.

A Lie $n$-algebra is \textbf{simple} if it is not one-dimensional and
every ideal $I\lhd V$ is either $0$ or $V$.

\begin{lemma}\label{le:simplequot}
  If $I\lhd V$ is a maximal ideal, then $V/I$ is simple or
  one-dimensional.
\end{lemma}

Simple Lie $n$-algebras have been classified.

\begin{theorem}[\cite{LingSimple}]\label{th:simple}
  A simple real Lie $n$-algebra is isomorphic to one of the
  ($n+1$)-dimensional Lie $n$-algebras defined, relative to a basis
  $\be_i$, by
  \begin{equation}
    \label{eq:simple-n-Lie}
    [\be_1,\dots, \widehat{\be_i},\dots, \be_{n+1}] = (-1)^i \varepsilon_i
    \be_i~,
  \end{equation}
  where a hat denotes omission and where the $\varepsilon_i$ are
  signs.
\end{theorem}

It is plain to see that simple real Lie $n$-algebras admit invariant
metrics of any signature.  Indeed, the Lie $n$-algebra in
\eqref{eq:simple-n-Lie} leaves invariant the diagonal metric with
entries $(\varepsilon_1, \dots, \varepsilon_{n+1})$.

The subspace $[V,\dots,V]$ is an ideal called the \textbf{derived
  ideal} of $V$.  Another ideal is provided by the \textbf{centre}
$Z$, defined by the condition $[Z,V,\dots, V]=0$.  More generally the
\textbf{centraliser} $Z(W)$ of a subspace $W \subset V$ is the
subalgebra defined by $[Z(W),W,V,\dots,V]=0$.

From now on let $(V,\Phi,b)$ be a metric Lie $n$-algebra.  If $W
\subset V$ is any subspace, we define
\begin{equation*}
  W^\perp = \left\{v \in V\middle | \left<v,w\right>=0~,\forall w\in
      W\right\}~.
\end{equation*}
Notice that $(W^\perp)^\perp = W$.  We say that $W$ is
\textbf{nondegenerate}, if $W \cap W^\perp = 0$, whence $V = W \oplus
W^\perp$; \textbf{isotropic}, if $W \subset W^\perp$; and
\textbf{coisotropic}, if $W \supset W^\perp$.  Of course, in
positive-definite signature, all subspaces are nondegenerate.

A metric Lie $n$-algebra is said to be \textbf{indecomposable} if it
is not isomorphic to a direct sum of orthogonal ideals or,
equivalently, if it does not possess any proper nondegenerate ideals:
for if $I\lhd V$ is nondegenerate, $V = I \oplus I^\perp$ is an
orthogonal direct sum of ideals.

\begin{lemma}\label{le:coisoquot}
  Let $I\lhd V$ be a coisotropic ideal of a metric Lie $n$-algebra.
  Then $I/I^\perp$ is a metric Lie $n$-algebra.
\end{lemma}

\begin{lemma}\label{le:centreperp}
  Let $V$ be a metric Lie $n$-algebra.  Then the centre is the
  orthogonal subspace to the derived ideal; that is,
  $[V,\dots,V]=Z^\perp$.
\end{lemma}

\begin{proposition}\label{pr:ideals}
  Let $V$ be a metric Lie $n$-algebra and $I \lhd V$ be an ideal.
  Then
  \begin{enumerate}
  \item $I^\perp \lhd V$ is also an ideal;
  \item $I^\perp\lhd Z(I)$; and
  \item if $I$ is minimal then $I^\perp$ is maximal.
  \end{enumerate}
\end{proposition}

\subsection{Structure of metric Lie $n$-algebras}
\label{sec:structure-metric-n}

We now investigate the structure of metric Lie $n$-algebras.  If a Lie
$n$-algebra is not simple or one-dimensional, then it has a proper
ideal and hence a minimal ideal.  Let $I\lhd V$ be a minimal ideal of
a metric Lie $n$-algebra.  Then $I \cap I^\perp$, being an ideal
contained in $I$, is either $0$ or $I$.  In other words, minimal
ideals are either nondegenerate or isotropic.  If nondegenerate, $V =
I \oplus I^\perp$ is decomposable.  Therefore if $V$ is
indecomposable, $I$ is isotropic.  Moreover, by
Proposition~\ref{pr:ideals} (2), $I$ is abelian and furthermore,
because $I$ is isotropic, $[I,I,V,\dots,V]=0$.

It follows that if $V$ is euclidean and indecomposable, it is either
one-dimensional or simple, whence of the form \eqref{eq:simple-n-Lie}
with all $\varepsilon_i=1$.  This result, originally due to
\cite{NagykLie}, was conjectured in \cite{FOPPluecker}.

Let $V$ be an indecomposable metric Lie $n$-algebra.  Then $V$ is
either simple, one-dimensional, or possesses an isotropic proper
minimal ideal $I$ which obeys $[I,I,V,\dots,V]=0$.  The perpendicular
ideal $I^\perp$ is maximal and hence by Lemma~\ref{le:simplequot}, $U
:= V/I^\perp$ is simple or one-dimensional, whereas by
Lemma~\ref{le:coisoquot}, $W:=I^\perp/I$ is a metric Lie $n$-algebra.

The inner product on $V$ induces a nondegenerate pairing $g: U \otimes
I \to \RR$.  Indeed, let $[u] = u + I^\perp \in U$ and $v\in I$.  Then
we define $g([u],v) = \left<u,v\right>$, which is clearly independent
of the coset representative for $[u]$.  In particular, $I \cong U^*$
is either one- or ($n+1$)-dimensional.  If the signature of the metric
of $W$ is $(p,q)$, that of $V$ is $(p+r,q+r)$ where $r = \dim I = \dim
U$.  So that if $V$ is to have lorentzian signature, $r=1$ and $W$
must be euclidean; although not necessarily indecomposable.

A lorentzian Lie $n$-algebra decomposes into one lorentzian
indecomposable factor and zero or more indecomposable euclidean
factors.  As discussed above, the indecomposable euclidean Lie
$n$-algebras are either one-dimensional or simple.  On the other hand,
an indecomposable lorentzian Lie $n$-algebra is either one-dimensional,
simple or else possesses a one-dimensional isotropic minimal ideal.
It is this latter case which remains to be treated and we do so now.

The quotient Lie $n$-algebra $U=V/I^\perp$ is also one-dimensional.
Let $u \in V$ be such that $u \not\in I^\perp$, whence its image in
$U$ generates it.  Because $I \cong U^*$, there is $v \in I$ such that
$\left<u,v\right> = 1$.  Complete it to a basis $(v,x_a)$ for
$I^\perp$.  Then $(u,v,x_a)$ is a basis for $V$, with $(x_a)$ spanning
a subspace isomorphic to $W=I^\perp/I$ and which, with a slight abuse
of notation, we will also denote $W$.  It is possible to choose $u$ so
that $\left<u,u\right> = 0$ and such that $\left<u,x\right>=0$ for all
$x\in W$.  Indeed, given any $u$, the map $x \mapsto \left<u,x\right>$
defines an element in the dual $W^*$.  Since the restriction of the
inner product to $W$ is nondegenerate, there is some $z\in W$ such
that $\left<u,x\right> = \left<z,x\right>$ for all $x \in W$.  We let
$u' = u - z$.  This still obeys $\left<u',v\right>=1$ and now also
$\left<u',x\right>=0$ for all $x\in W$.  Finally let $u'' = u' - \half
\left<u',u'\right> v$, which still satisfies $\left<u'',v\right> = 1$,
$\left<u'',x\right>=0$ for all $x\in W$, but now satisfies
$\left<u'',u''\right>=0$ as well.

From Proposition~\ref{pr:ideals} (2), it is immediate that
$[u,v,x_1\dots,x_{n-2}]=0=[v,x_1\dots,x_{n-1}]$, whence $v$ is
central.  Metricity then implies that the only nonzero $n$-brackets
take the form
\begin{equation}
  \label{eq:n-brackets}
  \begin{aligned}[m]
    [u,x_1,\dots, x_{n-1}] &=: [x_1,\dots,x_{n-1}]\\
    [x_1,\dots, x_n] &= (-1)^n \left<[x_1,\dots,x_{n-1}],x_n\right> v +
    [x_1,\dots, x_n]_W~.
  \end{aligned}
\end{equation}
The $n$-Jacobi identity is equivalent to the following two conditions:
\begin{enumerate}
\item $[x_1,\dots,x_{n-1}]$ defines a Lie ($n-1$)-algebra structure
  on $W$, which leaves the inner product invariant due to the
  skewsymmetry of $\left<[x_1,\dots,x_{n-1}],x_n\right>$; and
\item $[x_1,\dots, x_n]_W$ defines a euclidean Lie $n$-algebra
  structure on $W$ which is invariant under the ($n-1$)-algebra structure.
\end{enumerate}

We will show below that for $V$ indecomposable, the Lie $n$-algebra
structure on $W$ is abelian.

Let $(V,\Phi)$ be a Lie $n$-algebra.  It was already observed in
\cite{Filippov} that every $z\in V$ defines an $(n-1)$-bracket
$\Phi_z: \Lambda^{n-1}V \to V$, denoted simply as $[\dots]_z$, and
defined by
\begin{equation}
  \label{eq:n-1-Lie-z}
  [x_1,\dots,x_{n-1}]_z := [x_1,\dots,x_{n-1},z]~,
\end{equation}
which obeys the $(n-1)$-Jacobi identity as a consequence of the
$n$-Jacobi identity.  Thus $\Phi_z$ defines on $V$ a Lie
$(n-1)$-algebra structure for which $z$ is a central element.

If $(V,\Phi,b)$ is a metric Lie $n$-algebra, then each of the Lie
$(n-1)$-algebras $(V,\Phi_z,b)$ is a metric Lie $(n-1)$-algebra.  Let
$V$ be a simple euclidean Lie $n$-algebra given relative to a basis
$\be_i$ by \eqref{eq:simple-n-Lie} with all $\varepsilon_i =1$.
Moreover, such a basis is orthogonal, but not necessarily orthonormal.
Thus there is a one parameter family of such metric Lie $n$-algebras,
distinguished by the scale of the inner product.  We will denote the
simple Lie $n$-algebra with the above $n$-brackets by $\fs^{(n)}$.
Fixing any nonzero $x \in \fs^{(n)}$, the Lie $(n-1)$-algebra $\Phi_x$
is isomorphic to $\fs^{(n-1)} \oplus \RR$, where $\fs^{(n-1)}$ is the
simple euclidean Lie $(n-1)$-algebra structure on the perpendicular
complement of the line containing $x$.

\section{Lorentzian Lie $n$-algebras}
\label{sec:lorentzian}

We are now ready to classify the indecomposable lorentzian Lie
$n$-algebras and hence all lorentzian Lie $n$-algebras.  For $n=3$ the
classification of lorentzian Lie 3-algebras is given in \cite{Lor3Lie}
and is summarised in Theorem \ref{th:lorentzian3Lie}, whereas for
$n=2$, the lorentzian Lie algebras are classified in
\cite{MedinaLorentzian} and summarised in Theorem
\ref{th:lorentzianLie}.  Hence from now on we will take $n>3$.

We have previously shown that $V = \RR u \oplus \RR v \oplus W$, with
the $n$-brackets given by \eqref{eq:n-brackets}.  We will now show
that if $V$ is indecomposable, then as a Lie $n$-algebra, $W$ is
necessarily abelian.

Since $W$ is a euclidean Lie $n$-algebra, it is given by
\begin{equation*}
  W \cong \fa \oplus \fs^{(n)}_1 \oplus \dots \oplus \fs^{(n)}_q
\end{equation*}
where $\fa$ is a $p$-dimensional abelian Lie $n$-algebra.  The inner
product is such that the above direct sums are orthogonal, and the
inner product on each of the factors is positive-definite.

\begin{lemma}
  The Lie $(n-1)$-algebra structure on $W$ is such that the adjoint
  representation $\ad: \Lambda^{n-2}W \to \End W$ preserves the above
  orthogonal decomposition, whence each of the orthogonal summands are
  actually $(n-1)$-ideals.
\end{lemma}

\begin{proof}
  We will show that each of the $\fs^{(n)}_i$ are submodules.  Since
  the $(n-1)$-adjoint action of $W$ preserves the inner product, this
  means that so is $\fa$ and hence the claim, since submodules of the
  $(n-1)$-adjoint representation are $(n-1)$-ideals.

  Consider $\fs^{(n)}_1$, say, and let $(\be_1,\dots,\be_{n+1})$ be a
  basis.  For every $\bx := x_1 \wedge \dots \wedge x_{n-2} \in
  \Lambda^{n-2} W$,
  \begin{equation*}
    \ad_\bx \be_i := [x_1,\dots,x_{n-2},\be_i] = y_i + z_i~,
  \end{equation*}
  where $y_i \in \fs^{(n)}_1 \cap \be_i^\perp$ and $z_i \in
  \left(\fs^{(n)}_1\right)^\perp$.  We cannot have a component along
  $\be_i$ because invariance of the inner product says that $\ad_\bx
  \be_i \perp \be_i$.  Consider $\be_1 \wedge \dots \wedge
  \be_{n+1}$.  This being essentially the $n$-bracket in
  $\fs^{(n)}_1$, it is also preserved under $\ad_\bx$, whence
  \begin{align*}
    0 &= \ad_\bx \be_1 \wedge \dots \wedge \be_{n+1} + \dots + 
    \be_1 \wedge \dots \wedge \ad_\bx \be_{n+1}\\
    &= (y_1 + z_1) \wedge \dots \wedge \be_{n+1} + \dots + \be_1
    \wedge \dots \wedge (y_{n+1} + z_{n+1})\\
    &= z_1 \wedge \dots \wedge \be_{n+1} + \dots + \be_1 \wedge \dots
    \wedge z_{n+1}~.
  \end{align*}
  But each of these terms is independent, whence $z_i=0$ and $\ad_\bx$
  preserves $\fs^{(n)}_1$.
\end{proof}

As a Lie $(n-1)$-algebra, $W = W_0 \oplus W_1 \oplus \dots \oplus W_q$, where
$W_0$ is a $p$-dimensional euclidean Lie $(n-1)$-algebra and $W_{i>0}$
are $(n+1)$-dimensional euclidean Lie $(n-1)$-algebras.  Euclidean Lie
$(n-1)$-algebras are orthogonal direct sums of an abelian Lie
$(n-1)$-algebras (possibly zero-dimensional) and zero or more copies
of the $n$-dimensional euclidean simple Lie $(n-1)$-algebra
$\fs^{(n-1)}$.   In particular, on dimensional grounds, each $W_{i>0}$
is either abelian or isomorphic to $\fs^{(n-1)} \oplus \RR$.

We will now show that every $\fs^{(n)}$ summand in the Lie $n$-algebra
$W$ factorises in $V$, contradicting the assumption that $V$ is
indecomposable.

Consider one such $\fs^{(n)}$ summand, say $\fs^{(n)}_1$.  The
corresponding Lie $(n-1)$-algebra $W_1$ is either abelian or
isomorphic to $\fs^{(n-1)} \oplus \RR$.  If $W_1$ is abelian, so that
the $(n-1)$-brackets vanish, then it follows from
\eqref{eq:n-brackets} that for any $x \in \fs^{(n)}_1$,
$[u,x,V,\dots,V]=0$ and $[x,y,V,\dots,V]=0$ for any $y \in W$
perpendicular to $\fs^{(n)}_1$.  Hence $\fs^{(n)}_1\lhd V$ is a
nondegenerate ideal, contradicting the indecomposability of $V$.

If $W_1 \cong \fs^{(n-1)} \oplus \RR$, it has a one-dimensional centre
spanned by, say, $x\in W_1$.  Multiplying $x$ by a scalar if
necessary, we can assume that the Lie $(n-1)$-structure on
$W_1$ is given by the $(n-1)$-bracket $[x,\dots]_W$
induced from the Lie $n$-algebra structure on $\fs_1^{(n)}$; that is,
for all $y_1,\dots,y_{n-1}\in W_1$, $[y_1,\dots,y_{n-1}] =
[x,y_1,\dots,y_{n-1}]_W$.

This allows us to ``twist'' $\fs^{(n)}_1$ into a nondegenerate ideal
of $V$.   Indeed, define now a vector space isomorphism $\varphi: V
\to V$ by
\begin{equation}
  \label{eq:twist}
  \begin{aligned}[m]
    \varphi(v) &= v\\
    \varphi(u) &= u - x - \half |x|^2 v\\
    \varphi(y) &= y + \left<y,x\right> v\\
    \varphi(z) &= z~,
  \end{aligned}
\end{equation}
for all $y \in \fs^{(n)}_1$ and $z \in W\cap
\left(\fs^{(n)}_1\right)^\perp$.  It follows from
\eqref{eq:n-brackets} that, since $v$ is central, for all $y_i \in
\fs^{(n)}_1$,
\begin{align*}
  [\varphi(u),\varphi(y_1),\dots,\varphi(y_{n-1})] &= [u - x, y_1, \dots, y_{n-1}]\\
  &= [y_1,\dots,y_{n-1}] - [x,y_1,\dots,y_{n-1}]_W - (-1)^n \left<[x,y_1,\dots,y_{n-2}],y_{n-1}\right> v\\
  &= \left<[y_1,\dots,y_{n-2},y_{n-1}],x\right> v =0~,
\end{align*}
since $x$, being central, is perpendicular to the derived
$(n-1)$-ideal $[W_1,\dots,W_1]$ by Lemma \ref{le:centreperp}.
Finally, let $y_1,\dots,y_n \in \fs_1^{(n)}$, and since $v$ is
central, we have
\begin{align*}
  [\varphi(y_1),\dots,\varphi(y_n)] &= [y_1,\dots,y_n] \\
  &= (-1)^n \left<[y_1,\dots,y_{n-1}],y_n\right> v + [y_1,\dots,y_n]_W\\
  &= - \left<[y_1,\dots,y_{n-1},x]_W,y_n\right> v + [y_1,\dots,y_n]_W\\
  &=  \left<[y_1,\dots,y_{n-1},y_n]_W,x\right> v + [y_1,\dots,y_n]_W\\
  &= \varphi([y_1,\dots,y_n]_W)~.
\end{align*}
In other words, $\varphi(\fs^{(n)}_1)$ is a subalgebra.  Since it
commutes with $\varphi(u)$, $\varphi(v)$ and $\varphi(z)$ for $W \ni z
\perp \fs^{(n)}_1$, we see that $\varphi(\fs^{(n)}_1) \lhd V$.

It remains to show that $\varphi$ preserves the inner product.  If
$y\in\fs^{(n)}_1$, then
\begin{equation*}
  \left<\varphi(u), \varphi(y)\right> = \left<u - x- \half |x|^2 v, y +
    \left<x,y\right> v\right> = \left<x,y\right> \left<u,v\right> -
  \left<x,y\right> = 0~,
\end{equation*}
it is clear that $\left<\varphi(u),\varphi(v)\right> =1$ and that
$\left<\varphi(v),\varphi(v)\right> = 0$, whence we need only check
\begin{equation*}
  \left<\varphi(u),\varphi(u)\right> = \left<u-x-\half |x|^2 v, u-x-\half |x|^2
    v\right> = - |x|^2 \left<u,v\right> + \left<x,x\right> = 0~.
\end{equation*}
In other words, we conclude that $\varphi(\fs^{(n)}_1) \lhd V$ is a
degenerate ideal, contradicting again the fact that $V$ is
indecomposable.

Consequently there can be no $\fs^{(n)}$ summands in $W$, whence as a
Lie $n$-algebra, $W$ is abelian.  As a Lie $(n-1)$-algebra it is
euclidean, whence isomorphic to $\fb \oplus \fs^{(n-1)}_1 \oplus \dots
\oplus \fs^{(n-1)}_m$, where $\fb$ is now an abelian Lie
$(n-1)$-algebra.  However the abelian summand centralises $u$ and its
perpendicular complement in $W$, hence it is central and nondegenerate
in $V$, again contradicting the fact that it is indecomposable.
Therefore as a Lie $(n-1)$-algebra $W$ is euclidean semisimple: $W
\cong \fs^{(n-1)}_1 \oplus \dots \oplus \fs^{(n-1)}_m$.  In particular
it has dimension $mn$.

In summary, we have proved the following

\begin{theorem}\label{th:lorentzianNLie}
  Let $(V,\Phi,b)$ be a finite-dimensional indecomposable lorentzian
  Lie $n$-algebra.  Then it is either one-dimensional, simple, or else
  $V = \RR u \oplus \RR v \oplus W$, with $u,v$ complementary null
  directions perpendicular to $W$, such that the nonzero $n$-brackets
  take the form
  \begin{equation*}
    [u,x_1,\dots, x_{n-1}] = [x_1,\dots,x_{n-1}]
    \qquad\text{and}\qquad
    [x_1,\dots,x_n] = (-1)^n \left<[x_1,\dots,x_{n-1}],x_n\right> v~,
  \end{equation*}
  where $[x_1,\dots,x_{n-1}]$ makes $W$ into a euclidean Lie
  $(n-1)$-algebra isomorphic to $m$ copies of the simple euclidean Lie
  $(n-1)$-algebra $\fs^{(n-1)}$.  In particular, $\dim V = m n + 2$.
\end{theorem}

\begin{corollary}
  Let $(V,\Phi,b)$ be a lorentzian Lie $n$-algebra.  Then it is
  isomorphic to
  \begin{equation*}
    V = V_0 \oplus A \oplus \underbrace{\fs^{(n)} \oplus\dots \oplus \fs^{(n)}}_q~,
  \end{equation*}
  where $A$ is a ($p\geq0$)-dimensional euclidean abelian Lie
  $n$-algebra, $q\geq 0$, and $V_0$ is either one-dimensional with a
  negative-definite inner product, a lorentzian simple Lie
  $n$-algebra, or an indecomposable lorentzian Lie $n$-algebra with
  brackets given in Theorem~\ref{th:lorentzianNLie}.
\end{corollary}

It will not have escaped the reader's attention that Theorems
\ref{th:lorentzianLie}, \ref{th:lorentzian3Lie} and
\ref{th:lorentzianNLie} are very similar.  In all cases an
indecomposable lorentzian Lie $n$-algebra is either one-dimensional,
isomorphic to a unique simple Lie $n$-algebra or else it is obtained
by a ``double extension'' from a euclidean semisimple Lie
$(n-1)$-algebra, where the rôle of a semisimple Lie $1$-algebra is
played by an invertible skewsymmetric endomorphism.  Skew-diagonalising
the endomorphism, we see that it too is a direct sum of irreducibles.
Hence in all cases except for $n=2$, there is a unique simple
euclidean $n$-Lie algebra.  We will not formalise the notion of double
extension here, except to note that it clearly extends to Lie
$n$-algebras for $n>2$.  In the lorentzian examples, we double extend
by a one-dimensional Lie $n$-algebra, but it is possible to double
extend by other Lie $n$-algebras as well.  For example, let
$\fs$ be a simple Lie $n$-algebra and consider the vector space $\fs
\oplus \fs^*$ with the following $n$-bracket.  If $x_i\in \fs$ and
$\alpha \in \fs^*$, we define $[x_1,\dots,x_n] =
[x_1,\dots,x_n]_{\fs}$ to be the $n$-bracket in $\fs$ and
$[x_1,\dots,x_{n-1},\alpha] =: \beta \in \fs^*$ be defined by
$\beta(x_n) = - \alpha([x_1,\dots,x_n]_{\fs})$.  This is a metric Lie
$n$-algebra of signature $(n+1,n+1)$, relative to the inner product
consisting of any invariant inner product on $\fs$, the dual pairing
between $\fs$ and $\fs^*$ and declared to be identically zero on
$\fs^*$.

\bibliographystyle{utphys}
\bibliography{AdS,AdS3,ESYM,Sugra,Geometry,Algebra}

\end{document}